\documentclass{amsart}
\usepackage{amssymb,latexsym}
\usepackage{amsfonts}

\newtheorem{thm}{Theorem}[section]
\newtheorem{lem}[thm]{Lemma}

\newtheorem{rem}[thm]{Remark}
\newtheorem{pro}[thm]{Proposition}
\newtheorem{defi}[thm]{Definition}

\newcommand{\ba}{\begin{array}}
\newcommand{\ea}{\end{array}}

\def \qed{\cqfd}

\newcommand*{\QEDB}{\hfill\ensuremath{\square}}
\def\qed{\vbox{\hrule
\hbox{\vrule\hbox to 5pt{\vbox to 8pt{\vfil}\hfil}\vrule}\hrule}}

\newcommand{\beg}{\begin{eqnarray*}}
\newcommand{\begn}{\begin{eqnarray}}
\newcommand{\en}{\end{eqnarray*}}
\newcommand{\enn}{\end{eqnarray}}

\begin{document}
\title[Cusp K\"ahler-Ricci flow on compact K\"ahler manifold]{Cusp K\"ahler-Ricci flow on compact K\"ahler manifold}
\keywords{cusp K\"ahler-Ricci flow, conical K\"ahler-Ricci flow, cusp K\"ahler-Einstein metric.}
\author{Jiawei Liu}
\address{Jiawei Liu\\Beijing International Center for Mathematical Research\\ Peking University \\ Beijing 100871\\ China\\} \email{jwliu@math.pku.edu.cn}
\author{Xi Zhang}
\address{Xi Zhang\\Key Laboratory of Wu Wen-Tsun Mathematics\\ Chinese Academy of Sciences\\School of Mathematical Sciences\\
University of Science and Technology of China\\
Hefei, 230026, China\\ } \email{mathzx@ustc.edu.cn}
\thanks{AMS Mathematics Subject Classification. 53C55,\ 32W20.}
\thanks{}

\begin{abstract} In this paper, by limiting twisted conical K\"ahler-Ricci flows, we prove the long-time existence and uniqueness of cusp K\"ahler-Ricci flow on compact K\"ahler manifold $M$ which carries a smooth hypersurface $D$ such that the twisted canonical bundle $K_M+D$ is ample. Furthermore, we prove that this flow converge to a unique cusp K\"ahler-Einstein metric.

\end{abstract}

\maketitle
\section{Introduction}
\setcounter{equation}{0}

In this paper, we study a type of singular K\"ahler-Ricci flows which are obtained by limiting (twisted) conical K\"ahler-Ricci flows. Our motivation for considering the limit flows of conical flows is to study the existence of singular K\"ahler-Einstein metrics as the cone angles tend to $0$. In \cite{Tian0}, Tian anticipated that the complete Tian-Yau K\"ahler-Einstein metric on the complement of a divisor should be the limit of conical K\"ahler-Einstein metrics as the cone angles tend to $0$.

Let $M$ be a compact K\"ahler manifold with complex dimension $n$ and $D\subset M$ be a smooth hypersurface. Here, by supposing that the twisted canonical bundle $K_M+D$ is ample, we prove the long-time existence, uniqueness and convergence of cusp K\"ahler-Ricci flow by limiting twisted conical K\"ahler-Ricci flows. As an application, we show the existence of cusp K\"ahler-Einstein metric \cite{Kob84,TY87} by using cusp K\"ahler-Ricci flow.

The conical K\"ahler-Ricci flow was introduced to attack the existence problem of conical K\"ahler-Einstein metric. This equation was first proposed in Jeffres-Mazzeo-Rubinstein's paper (Section 2.5 in \cite{JMR}). Song-Wang made some conjectures on the relation between the convergence of conical K\"ahler-Ricci flow and the greatest Ricci lower bound of $M$ (conjecture 5.2 in \cite{SW}). The long-time existence, regularity and limit behaviour of conical K\"ahler-Ricci flow have been widely studied, see the works of Liu-Zhang \cite{JWLXZ,JWLXZ1}, Chen-Wang \cite{CW,CW1}, Wang \cite{YQW}, Shen \cite{LMSH1,LMSH2}, Edwards \cite{GEDWA}, Nomura \cite{Nomura},  Liu-Zhang \cite{JWLCJZ} and Zhang \cite{YSZ}.

By saying a closed positive $(1,1)$-current $\omega$ is conical K\"ahler metric with cone angle $2\pi \beta $ ($0<\beta\leq1$ ) along $D$, we mean that $D$ is locally given by $\{z^{n}=0\}$ and $\omega$ is asymptotically equivalent to model conical metric
\begin{eqnarray}
\sqrt{-1} \sum_{j=1}^{n-1} dz^{j}\wedge d\overline{z}^{j}+\frac{\sqrt{-1}dz^{n}\wedge d\overline{z}^{n}}{|z^{n}|^{2(1-\beta)}}.
\end{eqnarray}
And by saying a closed positive $(1,1)$-current $\omega$ is cusp K\"ahler metric along $D$, we mean that $D$ is locally given by $\{z^{n}=0\}$ and $\omega$ is asymptotically equivalent to model cusp metric
\begin{eqnarray}
\sqrt{-1} \sum_{j=1}^{n-1} dz^{j}\wedge d\overline{z}^{j}+\frac{\sqrt{-1}dz^{n}\wedge d\overline{z}^{n}}{|z^{n}|^{2}\log^2|z^{n}|^{2}}.
\end{eqnarray}

Let $\omega_0$ be a smooth K\"ahler metric on $M$ and satisfy $c_1(K_M)+c_1(D)=[\omega_0]$. We pick a section $s$ of $L_D$ cutting out this hypersurface, then we fix a smooth hermitian metric $h$ on $L_D$ and let $\theta$ be its curvature form. In \cite{JWLXZ1}, we proved the long-time existence, uniqueness, regularity and convergence of conical K\"ahler-Ricci flow with weak initial data $\omega_{\varphi_0}\in \mathcal{E}_{p}(M,\omega_{0})$ when $p>1$, where
 \begin{eqnarray*}
\mathcal{E}_{p}(M,\omega_{0})&=&\big\{\varphi\in\mathcal{E}(M,\omega_{0})\ |\ \frac{(\omega_{0}+\sqrt{-1}\partial\bar{\partial}\varphi)^{n}}{\omega_{0}^{n}}\in L^{p}(M,\omega_{0}^{n}) \big\},\\
\mathcal{E}(M,\omega_{0})&=&\big\{\varphi\in PSH(M,\omega_{0})\ |\ \int_{M} (\omega_{0}+\sqrt{-1}\partial\bar{\partial}\varphi)^{n}= \int_{M} \omega^{n}_{0}\big\}.
\end{eqnarray*}
Let $\alpha$ be a smooth closed $(1,1)$-form and $\hat{\omega}_\beta=\omega_0+\sqrt{-1}k\partial\bar{\partial}|s|_h^{2\beta}$. When $c_1(M)=\mu[\omega_0]+(1-\beta)c_1(D)+[\alpha]$ $(\mu\in\mathbb{R})$, by our arguments in \cite{JWLXZ1}, there exists a unique long-time solution of twisted conical K\"ahler-Ricci flow
\begin{eqnarray}\label{TCKRF}
\begin{cases}
  \frac{\partial \omega_{\beta}(t)}{\partial t}=-Ric(\omega_{\beta}(t))+\mu\omega_{\beta}(t)+(1-\beta)[D]+\alpha\\
  \\
  \omega_{\beta}(t)|_{t=0}=\omega_{\varphi_0}\\
  \end{cases}
  \end{eqnarray}
in the following sense:
\begin{itemize}
  \item  For any $[\delta, T]$ ($\delta , T>0$), there exists constant $C$ such that
   \begin{eqnarray*}
\ \ \ \ \ \ \ \ \ \ \ \ \ \ \ \ \ \ \ \ \ \ \ \ \ \ \ \ \ \ \ \ C^{-1}\hat{\omega}_\beta\leq\omega_{\beta}(t)\leq C\hat{\omega}_\beta\ \ \ \ \ \ \ \ \ \ \ \ on\ \ [\delta,T]\times (M\setminus D);
  \end{eqnarray*}
  \item  on $(0,\infty)\times(M\setminus D)$, $\omega_{\beta}(t)$ satisfies smooth twisted K\"ahler-Ricci flow;
  \item on $(0,\infty)\times M$, $\omega_{\beta}(t)$ satisfies equation (\ref{TCKRF}) in the sense of currents;
  \item there exists metric potential $\varphi_{\beta}(t)\in C^{0}\big([0,\infty)\times M\big)\cap C^{\infty}\big((0,\infty)\times (M\setminus D)\big)$ such that $\omega_{\beta}(t)=\omega_{0}+\sqrt{-1}\partial\bar{\partial}\varphi_{\beta}(t)$ and $\lim\limits_{t\rightarrow0^{+}}\|\varphi_{\beta}(t)-\varphi_{0}\|_{L^{\infty}(M)}=0$;
  \item on $[\delta, T]$, there exist constant $\alpha\in(0,1)$ and $C^{\ast}$ such that the above metric petential $\varphi_{\beta}(t)$ is $C^{\alpha}$ on $M$ with respect to $\omega_{0}$ and $\| \frac{\partial\varphi_{\beta}(t)}{\partial t}\|_{L^{\infty}(M\setminus D)}\leqslant C^{\ast}$.
  \end{itemize}
From Guenancia's result (Lemma $3.1$ in \cite{G11}),
\begin{eqnarray}
\omega_{\beta}=\omega_{0}-\sqrt{-1}\partial\bar{\partial}\log (\frac{1-|s|_{h}^{2\beta}}{\beta})^{2}:=\omega_{0}+\sqrt{-1}\partial\bar{\partial}\psi_{\beta}
\end{eqnarray}
is a conical K\"ahler metric  with cone angle $2\pi\beta$ along $D$. Hence, $\omega_\beta\in\mathcal{E}_{p}(M,\omega_{0})$ for $p\in(1,\frac{1}{1-\beta})$. By direct calculations, it is obvious that $\omega_\beta\geq\frac{1}{2}\omega_0$ for choosing suitable hermitian metric $h$ and $\omega_\beta$ converge to cusp K\"ahler metric
\begin{eqnarray}
\omega_{cusp}=\omega_{0}-\sqrt{-1}\partial\bar{\partial}\log\log^2|s|_h^2:=\omega_{0}+\sqrt{-1}\partial\bar{\partial}\psi_{0}
\end{eqnarray}
as $\beta\rightarrow 0$. After choosing $\mu=-1$, $\alpha=\beta\theta$ and $\omega_{\varphi_0}=\omega_\beta$ in $(\ref{TCKRF})$, we obtain twisted conical K\"ahler-Ricci flow
\begin{eqnarray}\label{TCKRF1}
\begin{cases}
  \frac{\partial \omega_{\beta}(t)}{\partial t}=-Ric(\omega_{\beta}(t))-\omega_{\beta}(t)+(1-\beta)[D]+\beta\theta.\\
  \\
  \omega_{\beta}(t)|_{t=0}=\omega_{\beta                                                                                                                                                                                           }\\
  \end{cases}
  \end{eqnarray}
Then by proving uniform estimates ( independent of $\beta$) for twisted conical K\"ahler-Ricci flows $(\ref{TCKRF1})$, we obtain a long-time solution to cusp K\"ahler-Ricci flow
\begin{eqnarray}\label{CUSPKRF}
\begin{cases}
  \frac{\partial \omega(t)}{\partial t}=-Ric(\omega(t))-\omega(t)+[D].\\
  \\
  \omega(t)|_{t=0}=\omega_{cusp}\\
  \end{cases}
  \end{eqnarray}

\begin{defi}\label{04.5}We call $\omega(t)$ a long-time solution to cusp K\"ahler-Ricci flow $(\ref{CUSPKRF})$ if it satisfies the following conditions.\\
  \textbf{(1)}  For any $[\delta, T]$ ($\delta , T>0$), there exists constant $C$ such that
   \begin{eqnarray*}
\ \ \ \ \ \ \ \ \ \ \ \ \ \ \ \ \ \ \ \ \ \ \ \ \ \ \ \ \ \ \ \ C^{-1}\omega_{cusp}\leq\omega(t)\leq C\omega_{cusp}\ \ \ \ \ \ \ \ \ \ \ \ on\ \ [\delta,T]\times(M\setminus D);
  \end{eqnarray*}
  \textbf{(2)} on $(0,\infty)\times(M\setminus D)$, $\omega(t)$ satisfies smooth K\"ahler-Ricci flow;\\
  \textbf{(3)} on $(0,\infty)\times M$, $\omega(t)$ satisfies equation (\ref{CUSPKRF}) in the sense of currents;\\
  \textbf{(4)} there exists $\varphi(t)\in C^{0}\big([0,\infty)\times (M\setminus D)\big)\cap C^{\infty}\big((0,\infty)\times (M\setminus D)\big)$ such that
 \begin{eqnarray*}\omega(t)=\omega_{0}+\sqrt{-1}\partial\bar{\partial}\varphi(t)\ \ and\ \ \lim\limits_{t\rightarrow0^{+}}\|\varphi(t)-\psi_{0}\|_{L^{1}(M)}=0; \end{eqnarray*}
  \textbf{(5)} on $(0,T]$, $\|\varphi(t)-\psi_0\|_{L^{\infty}(M\setminus D)}\leqslant C$;\\
  \textbf{(6)} on $[\delta, T]$, there exist constant $C$ such that $\| \frac{\partial\varphi(t)}{\partial t}\|_{L^{\infty}(M\setminus D)}\leqslant C$.
\end{defi}

 There are some important results on K\"ahler-Ricci flows (as well as its twisted versions with smooth twisting forms) from weak initial data, such as Chen-Ding\cite{CD}, Chen-Tian-Zhang \cite{CTZ}, Guedj-Zeriahi \cite{VGAZ2},  Lott-Zhang \cite{LZ,LZ1}, Nezza-Lu \cite{NL2014}, Song-Tian \cite{JSGT},  Sz\'ekelyhidi-Tosatti \cite{GSVT} and Zhang \cite{ZZ} etc. In particular, on $M\setminus D$, Lott-Zhang \cite{LZ} thoroughly studied the existence and convergence of K\"ahler-Ricci flow whose initial metric is finite volume K\"ahler metric with ``superstandard spatial asymptotics" ( Definition $8.10$ in \cite{LZ}). Their flow keeps ``superstandard spatial asymptotics" and this type of metrics contain cusp K\"ahler metrics. Here we consider K\"ahler-Ricci flow with non-smooth twisting form and weak initial data on $M$ globally, which can be seen as Lott-Zhang's case in some sense when we restrict it on $M\setminus D$. There are also some significant results on singular Ricci flows, see Ji-Mazzeo-Sesum \cite{JMS111}, Kleiner-Lott \cite{KLLO}, Mazzeo-Rubinstein-Sesum \cite{MRS}, Topping \cite{Top12} and Topping-Yin \cite{Top1010} etc.

 In \cite{JWLXZ1}, we studied conical K\"ahler-Ricci flow which is twisted by non-smooth twisting form and starts from weak initial data with $L^p$-density for $p>1$. Here, by approximating twisted conical K\"ahler-Ricci flows $(\ref{TCKRF1})$, we study cusp K\"ahler-Ricci flow with initial data $\omega_{cusp}$ which only admits $L^1$-density. For obtaining flow $(\ref{CUSPKRF})$, in addition to getting uniform estimates ( independent of $\beta$) of  flows $(\ref{TCKRF1})$, it is important to prove that $\varphi(t)$ converge to $\psi_{0}$ globally in $L^{1}$-sense and locally in $L^{\infty}$-sense outside $D$ as $t\rightarrow0^+$. In this process, we need to construct auxiliary function, and we also need a key observation (Proposition \ref{21888} and \ref{218}) that both $\psi_\beta$ and $\varphi_\beta(t)$ are monotone decreasing as $\beta\searrow0$. Then we obtain a uniqueness result of cusp K\"ahler-Ricci flow. In fact, we obtain the following theorem.

\begin{thm}\label{04} Let $M$ be a compact K\"ahler manifold and $\omega_{0}$ be a smooth K\"ahler metric. Assume that $D\subset M$ is a smooth hypersurface which satisfies $c_1(K_M)+c_1(D)=[\omega_0]$. Then there exists a unique long-time solution $\omega(t)=\omega_{0}+\sqrt{-1}\partial\bar{\partial}\varphi(t)$ to cusp K\"ahler-Ricci flow $(\ref{CUSPKRF})$.
\end{thm}

\begin{rem}\label{20161228} The uniqueness need to be understood in the following sense: if $\phi(t)\in C^{0}\big([0,\infty)\times (M\setminus D)\big)\cap C^{\infty}\big((0,\infty)\times (M\setminus D)\big)$ is a solution to equation
 \begin{eqnarray}\label{CMAE100}
\begin{cases}
 \frac{\partial \phi(t)}{\partial t}=\log\frac{(\omega_{0}+\sqrt{-1}\partial\bar{\partial}\phi(t))^{n}}{\omega_{0}^{n}}-\phi(t)+h_{0}+\log|s|_{h}^{2}\\
 \\
  \phi(0)=\psi_0\\
  \end{cases}
\end{eqnarray}
on $(0,\infty)\times (M\setminus D)$ and satisfies $(1)$, $(4)$, $(5)$ and $(6)$ in Definition $\ref{04.5}$, then $\phi(t)$ lies below $\varphi(t)$ which is obtained by limiting twisted conical K\"ahler-Ricci flows (\ref{TCKRF1}) in Theorem $\ref{04}$. When $n=1$, this uniqueness property is called ``maximally stretched" in Topping's \cite{Top10} and Giesen-Topping's \cite{GT11} works.
\end{rem}

\begin{rem} Since $K_M+D$ is ample, $K_{M}+(1-\beta)D$ is also ample for sufficiently small $\beta$. Guenancia \cite{G11} proved that cusp K\"ahler-Einstein metric is the limit of conical K\"ahler-Einstein metrics with background metrics $\omega_{0}-\beta\theta$ as $\beta\rightarrow0$. The cohomology classes are changing in this process. But in the flow case, we can not obtain a uniqueness result of cusp K\"ahler-Ricci flow $(\ref{CUSPKRF})$ if we choose the approximating flows that are conical  K\"ahler-Ricci flows with background metrics $\omega_{0}-\beta\theta$. In fact, if we choose the approximating flows that are conical K\"ahler-Ricci flows
\begin{eqnarray}\label{CKRF1}
\begin{cases}
  \frac{\partial \tilde{\omega}_{\beta}(t)}{\partial t}=-Ric(\tilde{\omega}_{\beta}(t))-\tilde{\omega}_{\beta}(t)+(1-\beta)[D]\\
  \\
  \tilde{\omega}_{\beta}(t)|_{t=0}=\omega_{0}-\beta\theta+\sqrt{-1}\partial\bar{\partial}\psi_\beta\\
  \end{cases}
  \end{eqnarray}
with background metrics $\omega_{0}-\beta\theta$, that is, $\tilde{\omega}_{\beta}(t)=\omega_{0}-\beta\theta+\sqrt{-1}\partial\bar{\partial}\tilde{\varphi}_{\beta}(t)$, we can also get a long-time solution $\tilde{\omega}(t)=\omega_{0}+\sqrt{-1}\partial\bar{\partial}\tilde{\varphi}(t)$ to equation $(\ref{CUSPKRF})$. But we do not know whether $\tilde{\varphi}(t)$ is unique or maximal. We can only prove $\tilde{\varphi}_\beta(t)+\beta\log|s|_h^2\nearrow\tilde{\varphi}(t)$ outside $D$ as $\beta\searrow0$. However, by the uniqueness result in Theorem $\ref{04}$,  $\tilde{\varphi}(t)$ must lie below $\varphi(t)$. Therefore, we set the background metric to $\omega_{0}$ in this paper.
\end{rem}

At last, we prove the convergence of cusp K\"ahler-Ricci flow $(\ref{CUSPKRF})$.
\begin{thm}\label{05} Cusp K\"ahler-Ricci flow $(\ref{CUSPKRF})$ converge to a K\"ahler-Einstein metric with cusp singularity along $D$ in $C_{loc}^{\infty}$-topology outside hypersurface $D$ and globally in the sense of currents.
\end{thm}

Kobayashi \cite{Kob84} and Tian-Yau \cite{TY87} asserted that if the twisted canonical bundle $K_M+D$ is ample, then there is a unique (up to constant multiple) complete cusp K\"ahler-Einstein metric with negative Ricci curvature on $M\setminus D$. The above convergence result recovers the existence of this cusp K\"ahler-Einstein metric.

The paper is organized as follows. In section $2$, we prove the long-time existence and uniqueness of cusp K\"ahler-Ricci flow $(\ref{CUSPKRF})$ by limiting twisted conical K\"ahler-Ricci flows (\ref{TCKRF1}) and constructing
auxiliary function. In section $3$, we prove the convergence theorem.
\medskip

{\bf  Acknowledgement:} The authors would like to thank Professor Jiayu Li for providing many suggestions and encouragement. The first author also would like to thank Professor Xiaohua Zhu  for his constant help. The second author is partially supported by NSF in China No.11625106, 11571332 and 11131007.

\section{The long-time existence of cusp K\"ahler-Ricci flow}
\setcounter{equation}{0}

In this section, we prove the long-time existence of cusp K\"ahler-Ricci flow by limiting twisted conical K\"ahler-Ricci flows (\ref{TCKRF1}), and we also prove the uniqueness theorem. For further consideration in the following arguments, we shall pay attention to the estimates which are independent of $\beta$.

From our arguments in \cite{JWLXZ1}, we know that there exists a unique long-time solution $\varphi_\beta(t)\in C^0\big([0,\infty)\times M\big)\bigcap C^\infty\big((0,\infty)\times(M\setminus D)\big)$ to the following conical parabolic complex Monge-Amp\`ere equation
\begin{eqnarray}\label{CMAE1}
\begin{cases}
 \frac{\partial \varphi_{\beta}(t)}{\partial t}=\log\frac{(\omega_{0}+\sqrt{-1}\partial\bar{\partial}\varphi_{\beta}(t))^{n}}{\omega_{0}^{n}}-\varphi_{\beta}(t)+h_{0}+(1-\beta)\log|s|_{h}^{2}\\
  \\
  \varphi_{\beta}(0)=\psi_\beta\\
  \end{cases}
\end{eqnarray}
on $(0,\infty)\times(M\setminus D)$, where $h_{0}$ satisfies $-Ric(\omega_{0})+\theta-\omega_{0}=\sqrt{-1}\partial\overline{\partial}h_{0}$. Let $\phi_{\beta}(t)=\varphi_{\beta}(t)-\psi_{\beta}$, we write the equation $(\ref{CMAE1})$ as
\begin{eqnarray}\label{CMAE2}
\begin{cases}
 \frac{\partial \phi_{\beta}(t)}{\partial t}=\log\frac{(\omega_{\beta}+\sqrt{-1}\partial\bar{\partial}\phi_{\beta}(t))^{n}}{\omega_{\beta}^{n}}-\phi_{\beta}(t)+h_{\beta}\\
  \phi_{\beta}(0)=0\\
  \end{cases}
\end{eqnarray}
on $(0,\infty)\times(M\setminus D)$, where $h_{\beta}=-\psi_{\beta}+h_{0}+\log\frac{|s|_{h}^{2(1-\beta)}\omega_{\beta}^{n}}{\omega_{0}^{n}}$ is uniformly bounded by constant $C$ independent of $\beta$.
\begin{lem}\label{2016002} There exists constant $C$ independent of $\beta$ and $t$ such that
\begin{eqnarray}\|\phi_{\beta}(t)\|_{L^{\infty}(M)}\leqslant C.
\end{eqnarray}
\end{lem}

{\bf Proof:}\ \ For any $\varepsilon>0$, we let $\chi_{\beta,\varepsilon}(t)=\phi_{\beta}(t)+\varepsilon\log|s|_{h}^{2}$. Since $\chi_{\beta,\varepsilon}(t)$ is smooth on $M\setminus D$, bounded from above and goes to $-\infty$ near D, it achieves its maximum on $M\setminus D$. Let $(t_{0}, x_{0})$ be the maximum point of $\chi_{\beta,\varepsilon}(t)$ on $[0,T]\times M$ with $x_{0}\in M\setminus D$. If $t_{0}=0$, then we have
\begin{eqnarray}\phi_{\beta}(t)\leqslant -\varepsilon\log|s|_{h}^{2}.
\end{eqnarray}
If $t_{0}\neq0$. At $(t_{0}, x_{0})$, we have
\begin{eqnarray*}0\leqslant\frac{\partial \chi_{\beta,\varepsilon}(t)}{\partial t}&=&\log\frac{(\omega_{\beta}+\sqrt{-1}\partial\bar{\partial}\phi_{\beta}(t))^{n}}{\omega_{\beta}^{n}}-\phi_{\beta}(t)+h_{\beta}\\
&=&\log\frac{(\omega_{\beta}+\sqrt{-1}\partial\bar{\partial}\chi_{\beta,\varepsilon}(t)+\varepsilon\theta)^{n}}{\omega_{\beta}^{n}}-\phi_{\beta}(t)+h_{\beta}\\
&\leq&n\log2-\phi_{\beta}(t)+C.
\end{eqnarray*}
Hence, $\phi_{\beta}(t_{0},x_{0})\leqslant C$ and
\begin{eqnarray}\phi_{\beta}(t)\leqslant C-\varepsilon\log|s|_{h}^{2},
\end{eqnarray}
where constant $C$ independent of $\beta$, $t$ and $\varepsilon$. Let $\varepsilon\rightarrow0$, we have
$\phi_{\beta}(t)\leqslant C$ on $M\setminus D$. Since $\phi_{\beta}(t)$ is continuous, $\phi_{\beta}(t)\leqslant C$ on $M$.

For the minimum, we can reproduce the same arguments with $\tilde{\chi}_{\beta,\varepsilon}(t)=\phi_{\beta}(t)-\varepsilon\log|s|_{h}^{2}$, and get $\phi_{\beta}(t)\geq C$ on $M$.\QEDB

\medskip

We now prove the uniform equivalence of volume forms along complex Monge-Amp\`ere equations $(\ref{CMAE2})$.

\begin{lem}\label{204} For any $T>0$, there exists constant $C$ independent of $\beta$ such that for any $t\in(0,T]$,
\begin{eqnarray}\label{201611303}\frac{t^n}{C}\leq\frac{(\omega_{\beta}+\sqrt{-1}\partial\bar{\partial}\phi_{\beta}(t))^{n}}{\omega_{\beta}^n}\leq e^{\frac{C}{t}}\ \ \ \ \ on\ \ M\setminus D.
\end{eqnarray}
\end{lem}

\medskip

{\bf Proof:}\ \ For any $t>0$, we assume that $t\in[\delta,T]$ with $\delta>0$. Let $\Delta_{\beta,t}$ be the Laplacian operator associated to  $\omega_\beta(t)=\omega_{\beta}+\sqrt{-1}\partial\bar{\partial}\phi_{\beta}(t)$. Straightforward calculations show that
\begin{eqnarray}(\frac{\partial}{\partial t}-\Delta_{\beta,t})\dot{\phi}_{\beta}(t)=-\dot{\phi}_{\beta}(t).
\end{eqnarray}

Let $H_{\beta,\varepsilon}^+(t)=(t-\delta)\dot{\phi}_{\beta}(t)-\phi_{\beta}(t)+\varepsilon\log|s|_h^2$. Since $H_{\beta,\varepsilon}^+(t)$ is smooth on $M\setminus D$, bounded from above and goes to $-\infty$ near D, it achieves its maximum on $M\setminus D$. Let $(t_{0}, x_{0})$ be the maximum point of $H_{\beta,\varepsilon}^+(t)$ on $[\delta,T]\times M$ with $x_{0}\in M\setminus D$. If $t_{0}=\delta$, then
\begin{eqnarray}(t-\delta)\dot{\phi}_{\beta}(t)\leq C-\varepsilon\log|s|_h^2,
\end{eqnarray}
where constant $C$ independent of $\beta$, $\delta$, $t$ and $\varepsilon$.
If $t_0\neq\delta$, then we have
\begin{eqnarray}\nonumber(\frac{\partial}{\partial t}-\Delta_{\beta,t})H_{\beta,\varepsilon}^+(t)&=&-(t-\delta)\dot{\phi}_{\beta}(t)+n+tr_{\omega_\beta(t)}(-\omega_\beta+\varepsilon\theta)\\
&\leqslant&-(t-\delta)\dot{\phi}_{\beta}(t)+n
\end{eqnarray}
for sufficiently small $\varepsilon$. By the maximum principle, we have
\begin{eqnarray}(t-\delta)\dot{\phi}_{\beta}(t)\leq C-\varepsilon\log|s|_h^2,
\end{eqnarray}
where constant $C$ independent of $\beta$, $\delta$, $t$ and $\varepsilon$. Let $\varepsilon\rightarrow0$ and then $\delta\rightarrow0$, we have
\begin{eqnarray}\label{201611301}\dot{\phi}_{\beta}(t)\leq \frac{C}{t}\ \ \ \ \ on\ \ (0,T]\times(M\setminus D),
\end{eqnarray}
where constant $C$ independent of $\beta$ and $t$.

Let $H_{\beta,\varepsilon}^-(t)=\dot{\phi}_{\beta}(t)+2\phi_{\beta}(t)-n\log (t-\delta)-\varepsilon\log|s|_h^2$. Then $H_{\beta,\varepsilon}^-(t)$ tend to $+\infty$ either $t\rightarrow0^+$ or $x\rightarrow D$.
\begin{eqnarray}(\frac{\partial}{\partial t}-\Delta_{\beta,t})H_{\beta,\varepsilon}^-(t)\geq\dot{\phi}_{\beta}(t)-2n-\frac{n}{t-\delta}+tr_{\omega_{\beta}(t)}\omega_\beta.
\end{eqnarray}
Assume that $(t_0,x_0)$ is the minimum point of $H_{\beta,\varepsilon}^-(t)$ on $[\delta,T]\times M$ with $t_0>\delta$ and $x_0\in M\setminus D$. There exists constant $C_1$ and $C_2$ such that
\begin{eqnarray}\label{2002}\nonumber(\frac{\partial}{\partial t}-\Delta_{\beta,t})H_{\beta,\varepsilon}^-(t)|_{(t_0,x_0)}&\geq& \big(C_1(\frac{\omega_{\beta}^n}{\omega^{n}_{\beta}(t)})^{\frac{1}{n}}+\log\frac{\omega^n_{\beta}(t)}{\omega_{\beta}^n}
-\frac{C_2}{t-\delta}\big)|_{(t_0,x_0)}\\
&\geq&\big(\frac{C_1}{2}(\frac{\omega_{\beta}^n}{\omega^{n}_{\beta}(t)})^{\frac{1}{n}}
-\frac{C_2}{t-\delta}\big)|_{(t_0,x_0)},
\end{eqnarray}
where constant $C_1$ depends only on $n$, $C_2$  depends only on $n$, $\omega_{0}$ and $T$. In inequality $(\ref{2002})$, without loss of generality, we assume that $\frac{\omega_{\beta}^n}{\omega^{n}_{\beta}(t)}>1$ and $\frac{C_1}{2}(\frac{\omega_{\beta}^n}{\omega^{n}_{\beta}(t)})^{\frac{1}{n}}+\log\frac{\omega^n_{\beta}(t)}
{\omega_{\beta}^n}\geq0$ at $(t_0,x_0)$. By the maximum principle, we have
\begin{eqnarray}\omega^{n}_{\beta}(t_0,x_0)\geq C_4(t_0-\delta)^n \omega_{\beta}^n(x_0),
\end{eqnarray}
where $C_4$ independent of $\beta$, $\varepsilon$ and $\delta$. Then it easily follows that
\begin{eqnarray}\dot{\phi}_\beta(t)\geq-C+n\log(t-\delta)+\varepsilon\log|s|_h^2,
\end{eqnarray}
where constant $C$ independent of $\beta$, $\varepsilon$ and $\delta$. Let $\varepsilon\rightarrow0$ and then $\delta\rightarrow0$, we have
\begin{eqnarray}\label{201611302}\dot{\phi}_{\beta}(t)\geq -C+n\log t\ \ \ \ \ on\ \ (0,T]\times(M\setminus D),
\end{eqnarray}
where constant $C$ independent of $\beta$. By $(\ref{201611301})$ and $(\ref{201611302})$, we obtain $(\ref{201611303})$. \QEDB

\medskip

We first recall Guenancia's results about the curvature of $\omega_\beta$ ( Theorem $3.2$ \cite{G11}).
\begin{lem}\label{2016001} There exists a constant $C$ depending only on $M$ such that for all $\beta\in(0,\frac{1}{2}]$, the holomorphic bisectional curvature of $\omega_{\beta}$ is bounded by $C$.
\end{lem}

Next, we prove the uniform equivalence of metrics along twisted conical K\"ahler-Ricci flows $(\ref{TCKRF1})$ by Chern-Lu inequality.

\begin{lem}\label{205} For any $T>0$, there exists constant $C$ independent of $\beta$ such that for any $t\in(0,T]$,
\begin{eqnarray}e^{-\frac{C}{t}}\omega_\beta\leq\omega_{\beta}(t)\leq e^{\frac{C}{t}}\omega_\beta\ \ \ \ \ on\ \ M\setminus D.
\end{eqnarray}
\end{lem}

{\bf Proof:}\ \ By Chern-Lu inequality \cite{SSC, YCL}, on $M\setminus D$, we have
\begin{eqnarray}\nonumber\Delta_{\beta,t}\log tr_{\omega_{\beta}(t)}\omega_{\beta}&=&\frac{g_{\beta t}^{i\bar{l}}g_{\beta t}^{k\bar{j}}R_{\beta ti\bar{j}}g_{\beta k\bar{l}}}{tr_{\omega_{\beta}(t)}\omega_{\beta}}
-\frac{1}{tr_{\omega_{\beta}(t)}\omega_{\beta}}(g_{\beta t}^{i\bar{j}}g_{\beta t}^{k\bar{l}}R_{\beta i\bar{j}k\bar{l}})\\
&-&\frac{g_{\beta t}^{i\bar{j}}\partial_{i}tr_{\omega_{\beta}(t)}\omega_{\beta}\partial_{\bar{j}}tr_{\omega_{\beta}(t)}\omega_{\beta} }{(tr_{\omega_{\beta}(t)}\omega_{\beta})^{2}}+\frac{g_{\beta t}^{i\bar{j}}g_{\beta t}^{k\bar{l}}g_{\beta}^{p\bar{q}}\partial_{i}g_{\beta k\bar{q}}\partial_{\bar{j}}g_{\beta p\bar{l}}}{tr_{\omega_{\beta}(t)}\omega_{\beta}}.
\end{eqnarray}
At the same time, on $M\setminus D$,
\begin{eqnarray}\frac{\partial}{\partial t}\log tr_{\omega_{\beta}(t)}\omega_{\beta}=\frac{g_{\beta t}^{i\bar{l}}g_{\beta t}^{k\bar{j}}g_{\beta i\bar{j}}(R_{\beta tk\bar{l}}+g_{\beta tk\bar{l}}-\beta\theta_{k\bar{l}})}{tr_{\omega_{\beta}(t)}\omega_{\beta}}.
\end{eqnarray}
By using Cauchy-Schwarz inequality and Lemma \ref{2016001}, we have
\begin{eqnarray}\nonumber&\ &(\frac{\partial}{\partial t}-\Delta_{\beta,t})\log tr_{\omega_{\beta}(t)}\omega_{\beta}\leq C tr_{\omega_{\beta}(t)}\omega_{\beta}+1,
\end{eqnarray}
where constant $C$ independent of $\beta$.

Let $H_{\beta,\varepsilon}(t)=(t-\delta)\log tr_{\omega_{\beta}(t)}\omega_{\beta}-A\phi_{\beta}(t)+\varepsilon\log|s|_h^2$,  $A$ be a sufficiently large constant and $(t_0,x_0)$ be the maximum point of $H_{\beta,\varepsilon}(t)$ on $[\delta,T]\times (M\setminus D)$. We know that $x_0 \in M\setminus D$ and we need only consider $t_0>\delta$. By direct calculations,
\begin{eqnarray}\nonumber(\frac{\partial}{\partial t}-\Delta_{\beta,t})H_{\beta,\varepsilon}(t)&\leq&\log tr_{\omega_{\beta}(t)}\omega_{\beta}+Ctr_{\omega_{\beta}(t)}\omega_{\beta}-A\dot{\phi}_{\varepsilon,j}(t)-Atr_{\omega_{\beta}(t)}\omega_{\beta}+\varepsilon tr_{\omega_{\beta}(t)}\theta+C\\\nonumber
&\leq&-\frac{A}{2}tr_{\omega_{\beta}(t)}\omega_{\beta}+\log tr_{\omega_{\beta}(t)}\omega_{\beta}-A\log\frac{\omega_{\beta}^n(t)}{\omega_{\beta}^n}+C,
\end{eqnarray}
where constant $C$ independent of $\beta$ and $\delta$.

Without loss of generality, we assume that
$-\frac{A}{4}tr_{\omega_{\beta}(t)}\omega_{\beta}+\log tr_{\omega_{\beta}(t)}\omega_{\beta}\leq0$ at $(t_0,x_0)$. Then at $(t_0,x_0)$, by Lemma \ref{204}, we have
\begin{eqnarray}(\frac{\partial}{\partial t}-\Delta_{\beta,t})H_{\beta,\varepsilon}(t)\leq-\frac{A}{4}tr_{\omega_{\beta}(t)}\omega_{\beta}-An\log (t-\delta)+C.
\end{eqnarray}
By the maximum principle, at $(t_0,x_0)$,
\begin{eqnarray}tr_{\omega_{\beta}(t)}\omega_{\beta}\leq C\log \frac{1}{t-\delta}+C,
\end{eqnarray}
which implies that
\begin{eqnarray}(t-\delta)\log tr_{\omega_{\beta}(t)}\omega_{\beta}\leq (t_0-\delta)\log(C\log \frac{1}{t_0-\delta}+C)+C-\varepsilon\log|s|_h^2.
\end{eqnarray}
Let $\varepsilon\rightarrow0$ and then $\delta\rightarrow0$, on $(0,T]\times (M\setminus D)$,
\begin{eqnarray}\label{208} tr_{\omega_{\beta}(t)}\omega_{\beta}\leq e^\frac{C}{t}.
\end{eqnarray}
By using inequality
\begin{eqnarray}tr_{\omega_{\beta}}\omega_{\beta}(t)\leq\frac{1}{(n-1)!}(tr_{\omega_{\beta}(t)}\omega_{\beta})^{n-1}\frac{\omega_{\beta}^{n}(t)}{\omega_{\beta}^{n}},
\end{eqnarray}
we have
\begin{eqnarray}\label{207}tr_{\omega_{\beta}}\omega_{\beta}(t)\leq e^{\frac{C}{t}},
\end{eqnarray}
where $C$ independent of $\beta$. From $(\ref{208})$ and $(\ref{207})$, we prove the lemma.\QEDB

 By the argument as that in \cite{JWLXZ1}, we get the following local Calabi's $C^3$-estimates and curvature estimates.

\begin{lem} For any $T>0$ and $B_r(p)\subset\subset M\setminus D$, there exist constants $C$, $C'$ and $C''$ depend only on $n$, $T$, $\omega_0$ and $dist_{\omega_0}(B_r(p),D)$ such that
\begin{eqnarray*}S_{\omega_{\beta}(t)}&\leq&\frac{C'}{r^{2}}e^{\frac{C}{t}},\\
|Rm_{\omega_{\beta}(t)}|_{\omega_{\beta}(t)}^{2}&\leq&\frac{C''}{r^{4}}e^{\frac{C}{t}}
\end{eqnarray*}
on $(0,T]\times B_{\frac{r}{2}}(p)$.\end{lem}

By using the standard parabolic Schauder regularity theory \cite{GLIE}, we obtain the following proposition.

\begin{pro}\label{217} For any $0<\delta<T<\infty$, $k\in\mathbb{N}^{+}$ and $B_r(p)\subset\subset M\setminus D$, there exists constant $C_{\delta,T,k,p,r}$ depends only on $n$, $\delta$, $k$, $T$, $\omega_0$ and $dist_{\omega_{0}}(B_r(p),D)$ such that for $\beta\in(0,\frac{1}{2}]$,
\begin{eqnarray}\|\varphi_{\beta}(t)\|_{C^{k}\big([\delta,T]\times B_r(p)\big)}\leq C_{\delta,T,k,p,r}.
\end{eqnarray}
\end{pro}

Through a further observation to $\psi_\beta$ and equation $(\ref{CMAE1})$, we prove the monotonicity of $\psi_\beta$ and $\varphi_\beta(t)$ with respect to $\beta$.
\begin{pro}\label{21888} For any $x\in M$, $\psi_\beta(x)$ is monotone decreasing as $\beta\searrow0$.
\end{pro}

{\bf Proof:}\ \ By direct computations, for any $x\in M\setminus D$, we have
\begin{eqnarray}\frac{d\psi_\beta}{d\beta}=2\frac{\beta|s|_h^{2\beta}\log|s|_h^2+1-|s|_h^{2\beta}}{\beta(1-|s|_h^{2\beta})}.
\end{eqnarray}
Denote $f_\beta(x)=\beta x^\beta \log x+1-x^\beta$ for $\beta>0$ and $x\in[0,1]$.
\begin{eqnarray}f'_{\beta}(x)=\beta^2x^{\beta-1}\log x\leq 0.
\end{eqnarray}
Hence $f_\beta(x)\geq f_\beta(1)=0$ and we have $\frac{d\psi_\beta}{d\beta}\geq0$. \QEDB

\begin{pro}\label{218} For any $(t,x)\in (0,\infty)\times M$, $\varphi_\beta(t,x)$ is monotone decreasing as $\beta\searrow0$.
\end{pro}

{\bf Proof:}\ \ By the arguments in \cite{JWLXZ1}, we obtain $(\ref{CMAE1})$ by approximating equations
\begin{eqnarray}\label{CMAE101}
\begin{cases}
 \frac{\partial \varphi_{\beta\varepsilon}(t)}{\partial t}=\log\frac{(\omega_{0}+\sqrt{-1}\partial\bar{\partial}\varphi_{\beta\varepsilon}(t))^{n}}{\omega_{0}^{n}}-\varphi_{\beta\varepsilon}(t)+h_{0}+\log(\varepsilon^2+|s|_{h}^{2})^{1-\beta}\\
  \\
  \varphi_{\beta\varepsilon}(0)=\psi_\beta\\
  \end{cases}
\end{eqnarray}
For $\beta_1<\beta_2$, let $\psi_{1,2}(t)=\varphi_{\beta_1\varepsilon}(t)-\varphi_{\beta_2\varepsilon}(t)$. On $[\eta,T]\times M$ with $\eta>0$ and $T<\infty$,
\begin{eqnarray}\nonumber&\ &\frac{\partial }{\partial t}(e^{t-\eta}\psi_{1,2}(t))\\
&\leq&e^{t-\eta}\log\frac{\big(e^{t-\eta}\omega_{0}+\sqrt{-1}\partial\bar{\partial}e^{t-\eta}\varphi_{\beta_2\varepsilon}(t)+\sqrt{-1}\partial\bar{\partial}e^{t-\eta}\psi_{1,2}(t)\big)^{n}}
  {(e^{t-\eta}\omega_{0}+\sqrt{-1}\partial\bar{\partial}e^{t-\eta}\varphi_{\beta_2\varepsilon}(t))^{n}}.
\end{eqnarray}

Let $\tilde{\psi}_{1,2}(t)=e^{t-\eta}\psi_{1,2}(t)-\delta (t-\eta)$ with $\delta>0$ and $(t_0, x_0)$ be the maximum point of $\tilde{\psi}_{1,2}(t)$ on $[\eta,T]\times M$. If $t_0>\eta$, by the maximum principle, at this point,
\begin{eqnarray}0\leq\frac{\partial }{\partial t}\tilde{\psi}_{1,2}(t)=\frac{\partial }{\partial t}(e^{t-\eta}\psi_{1,2}(t))-\delta\leq-\delta
\end{eqnarray}
which is impossible, hence $t_0=\eta$. So for any $(t,x)\in [\eta,T]\times M$,
\begin{eqnarray}\psi_{1,2}(t,x)\leq e^{-t+\eta}\sup\limits_{M}\psi_{1,2}(\eta,x)+T\delta.
\end{eqnarray}
Since $\lim\limits_{t\rightarrow0^+}\parallel\varphi_{\beta\varepsilon}(t)-\chi_\beta\parallel_{L^\infty(M)}=0$, let $\eta\rightarrow0$, we get
\begin{eqnarray}\psi_{1,2}(t,x)\leq e^{-t}\sup\limits_{M}(\psi_{\beta_1}-\psi_{\beta_2})+T\delta\leq T\delta.
\end{eqnarray}
Let $\delta\rightarrow0$ and then $\varepsilon\rightarrow0$, we conclude that $\varphi_{\beta_1}(t,x)\leq\varphi_{\beta_2}(t,x)$.\QEDB

\medskip

For any $[\delta, T]\times K\subset\subset (0,\infty)\times M\setminus D$ and $k\geq0$, $\|\varphi_{\beta}(t)\|_{C^{k}([\delta,T]\times K)}$ is uniformly bounded by Proposition \ref{217}. Let $\delta$ approximate to $0$, $T$ approximate to $\infty$ and $K$ approximate to $M\setminus D$, by diagonal rule, we get a sequence $\{\beta_i\}$, such that $\varphi_{\beta_i}(t)$ converge in $C^\infty_{loc}$-topology on $(0,\infty)\times (M\setminus D)$ to a function $\varphi(t)$ that is smooth on $C^\infty\big((0,\infty)\times(M\setminus D)\big)$ and satisfies equation
\begin{eqnarray}\label{223}  \frac{\partial \varphi(t)}{\partial t}=\log\frac{(\omega_{0}+\sqrt{-1}\partial\bar{\partial}\varphi(t))^{n}}{\omega_{0}^{n}}-\varphi(t)+h_0
  +\log|s|_{h}^{2}
\end{eqnarray}
on $(0,\infty)\times (M\setminus D)$. Since $\varphi_\beta(t)$ is monotone decreasing as $\beta\rightarrow0$, $\varphi_{\beta}(t)$ converge in $C^\infty_{loc}$-topology on $(0,\infty)\times (M\setminus D)$ to $\varphi(t)$. For any $T>0$,
\begin{eqnarray}
\label{222}&\ &e^{-\frac{C}{t}}\omega_{cusp}\leq\omega(t)\leq e^{\frac{C}{t}}\omega_{cusp}\ \ on\ (0,T]\times (M\setminus D),
\end{eqnarray}
where $\omega(t)=\omega_0+\sqrt{-1}\partial\overline{\partial}\varphi(t)$, constants $C$ depend only on $n$, $\omega_{0}$ and $T$.

\medskip

Next, by using the monotonicity of $\varphi_{\beta}(t)$ with respect to $\beta$ and constructing auxiliary function, we prove the $L^1$-convergence of $\varphi(t)$ as $t\rightarrow0^{+}$ as well as $\varphi(t)$ converge to $\psi_{0}$ in $L^{\infty}$-norm as $t\rightarrow0^{+}$ on any compact subset $K\subset\subset M\setminus D$.

\begin{lem}\label{101010} There exists a unique $\varphi_\beta\in PSH(M,\omega_0)\bigcap L^\infty(M)$ to equation
\begin{eqnarray}\label{1020101}(\omega_0+\sqrt{-1}\partial\bar{\partial}\varphi_\beta)^n=e^{\varphi_\beta-h_0}\frac{\omega_0^n}{|s|_h^{2(1-\beta)}}.
\end{eqnarray}
Furthermore, $\varphi_\beta\in C^{2,\alpha,\beta}(M)$ and $\parallel\varphi_\beta-\psi_\beta\parallel_{L^\infty(M)}$ can be uniformly bounded by constant $C$ independent of $\beta$.
\end{lem}

{\bf Proof:}\ \ By Eyssidieux-Guedj-Zeriahi's theorem (see Theorem $4.1$ in \cite{EGZ}), there exists a unique continuous solution $\varphi_\beta$ to equation $(\ref{1020101})$. Then by Guenancia-P$\breve{a}$un's regularity estimates \cite{GP1} ( see also Liu-Zhang \cite{JWLCJZ}), $\varphi_\beta\in C^{2,\alpha,\beta}(M)$. Next, we prove $\parallel\varphi_\beta-\psi_\beta\parallel_{L^\infty(M)}$ can be uniformly bounded. Let $u_\beta=\varphi_\beta-\psi_\beta$, we write equation $(\ref{1020101})$ as
\begin{eqnarray}\label{1020102}(\omega_\beta+\sqrt{-1}\partial\bar{\partial}u_\beta)^n=e^{u_\beta+h_\beta}\omega_\beta^n,
\end{eqnarray}
where $h_{\beta}=\psi_{\beta}-h_{0}+\log\frac{\omega_{0}^{n}}{|s|_{h}^{2(1-\beta)}\omega_{\beta}^{n}}$ is uniformly bounded independent of $\beta$. Define $\chi_{\beta,\varepsilon}=u_\beta+\varepsilon\log|s|_h^2$. Then $\sqrt{-1}\partial\bar{\partial}\chi_{\beta,\varepsilon}=\sqrt{-1}\partial\bar{\partial}u_\beta-\varepsilon\theta$. Since $\chi_{\beta,\varepsilon}$ is smooth on $M\setminus D$, bounded from above and goes to $-\infty$ near D, it achieves its maximum on $M\setminus D$. Let $x_{0}$ be the maximum point of $\chi_{\beta,\varepsilon}$ on $M$ with $x_{0}\in M\setminus D$.
\begin{eqnarray}(\omega_\beta+\sqrt{-1}\partial\bar{\partial}u_\beta)^n(x_0)=(\omega_\beta+\sqrt{-1}\partial\bar{\partial}\chi_{\beta,\varepsilon}+\varepsilon\theta)^n(x_0)\leq2^n\omega_\beta^n(x_0).
\end{eqnarray}
By the maximum principle, $u_\beta\leq C-\varepsilon\log|s|_h^2$, where constant $C$ independent of $\beta$ and $\varepsilon$. Let $\varepsilon\rightarrow0$, we get the uniform upper bound of $u_\beta$. By the similar arguments, we can obtain the uniform lower bound of $u_\beta$.\QEDB

\begin{pro}\label{101} $\varphi(t)\in C^0([0,\infty)\times (M\setminus D))$ and
\begin{eqnarray}\label{102}\lim\limits_{t\rightarrow0^+}\|\varphi(t)-\psi_{0}\|_{L^{1}(M)}=0.
\end{eqnarray}
\end{pro}

{\bf Proof:}\ \ By the monotonicity of $\varphi_{\beta}(t)$ with respect to $\beta$, for any $(t,z)\in(0,T]\times (M\setminus D)$, we have
\begin{eqnarray}\label{2210}\nonumber\varphi(t,z)-\psi_0(z)&\leq&\varphi_{\beta}(t,z)-\psi_0(z)\\
&\leq&|\varphi_{\beta}(t,z)-\psi_{\beta}(z)|+|\psi_{\beta}(z)-\psi_{0}(z)|.
\end{eqnarray}
Since $\psi_\beta$ converge to $\psi_0$ in $C^\infty_{loc}$-sense outside $D$ as $\beta\rightarrow0$, and
\begin{eqnarray}\label{2016120201}\lim\limits_{t\rightarrow0^+}\|\varphi_{\beta}(t,z)-\psi_\beta\|_{L^{\infty}(M)}=0.
\end{eqnarray}
For any $\epsilon>0$ and $K\subset\subset M\setminus D$, there exists $N$ such that for $\beta_1<\frac{1}{N}$,
\begin{eqnarray}
\|\psi_{\beta_1}(z)-\psi_{0}(z)\|_{L^{\infty}(K)}&<&\frac{\epsilon}{2}.
\end{eqnarray}
Fix such $\beta_{1}$, there exists $0<\delta_1<T$ such that
\begin{eqnarray}\sup\limits_{[0,\delta_1]\times M}|\varphi_{\beta_1}(t,z)-\psi_{\beta_1}|<\frac{\epsilon}{2}.
\end{eqnarray}
Combining the above estimates together, for any $t\in(0,\delta_1]$ and $z\in K$
\begin{eqnarray}\sup\limits_{[0,\delta_1]\times K}(\varphi(t,z)-\psi_0(z))<\epsilon.
\end{eqnarray}

We define function
\begin{eqnarray}H_{\beta}(t)=(1-te^{-t})\psi_\beta+te^{-t}\varphi_{\beta}+h(t)e^{-t},
\end{eqnarray}
where $\varphi_\beta$ and $u_\beta=\varphi_\beta-\psi_\beta$ are obtained in Lemma $\ref{101010}$, and
\begin{eqnarray}\nonumber h(t)=(1-e^t-t)\|u_{\beta}\|_{L^\infty(M)}+n(t\log t-t)e^{t}-n\int_0^te^{s}s\log s ds.
\end{eqnarray}
Straightforward calculations show that
\begin{eqnarray*}\frac{\partial}{\partial t}H_{\beta}(t)+H_{\beta}(t)&=&\psi_\beta+e^{-t}u_\beta-e^{-t}\|u_{\beta}\|_{L^\infty(M)}-\|u_{\beta}\|_{L^\infty(M)}+n\log t-nt\\
&\leq&\psi_\beta+u_\beta+n\log t-nt\\
&\leq&\varphi_\beta+n\log t-nt
\end{eqnarray*}
Therefore, we have
\begin{eqnarray*}e^{\frac{\partial}{\partial t}H_{\beta}(t)+H_{\beta}(t)}\omega_0^n\leq
t^ne^{-nt}e^{\varphi_{\beta}}\omega_0^n.
\end{eqnarray*}
When $t$ is sufficiently small,
\begin{eqnarray*}\omega_0+\sqrt{-1}\partial\overline{\partial}H_{\beta}(t)&=&(1-te^{-t})(\omega_0+\sqrt{-1}\partial\overline{\partial}\psi_\beta)+te^{-t}(\omega_0+\sqrt{-1}\partial\overline{\partial}\varphi_{\beta})\\
&\geq&te^{-t}(\omega_0+\sqrt{-1}\partial\overline{\partial}\varphi_\beta).
\end{eqnarray*}
Combining the above inequalities,
\begin{eqnarray*}(\omega_0+\sqrt{-1}\partial\overline{\partial}H_{\beta}(t))^n&\geq&t^ne^{-n t}(\omega_0+\sqrt{-1}\partial\overline{\partial}\varphi_\beta)^n\\
&\geq&e^{-h_{0}+\frac{\partial}{\partial t}H_{\beta}(t)+H_{\beta}(t)}\frac{\omega_{0}^{n}}
{|s|_{h}^{2(1-\beta)}},
\end{eqnarray*}
which is equivalent to
\begin{eqnarray}
\begin{cases}
  \frac{\partial}{\partial t}H_{\beta}(t)\leq \log \frac{(\omega_0+\sqrt{-1}\partial\overline{\partial}H_{\beta}(t))^n}{\omega_{0}^{n}}-H_{\beta}(t)+h_{0}+\log|s|_{h}^{2(1-\beta)}.\\
  \\
  H_{\beta}(0)=\psi_\beta
  \end{cases}
\end{eqnarray}

Next, we prove $H_\beta(t)\leq \varphi_\beta(t)$ for sufficiently small $t$ by using Jeffres' trick \cite{TJEF}. For any $0<t_1<T<\infty$ and $a>0$.

Denote $\Psi(t)=H_\beta(t)+a|s|_h^{2q}-\varphi_\beta(t)$ and $\hat{\Delta}=\int_0^1 g_{sH_\beta(t)+(1-s)\varphi_\beta(t)}^{i\bar{j}}\frac{\partial^2}{\partial z^i\partial\bar{z}^j}ds$, where $0<q<1$ is determined later. $\Psi(t)$ evolves along the following equation
\begin{eqnarray*}
  \frac{\partial \Psi(t)}{\partial t}=\hat{\Delta}\Psi(t)-a\hat{\Delta}|s|_h^{2q}-\Psi(t)+a|s|_h^{2q}.
\end{eqnarray*}
Since
\begin{eqnarray*}\omega_0+\sqrt{-1}\partial\overline{\partial}H_{\beta}(t)&\geq&(1-te^{-t})(\omega_0+\sqrt{-1}\partial\overline{\partial}\psi_\beta)\geq \frac{1}{4}\omega_0,\\
\omega_0+\sqrt{-1}\partial\overline{\partial}\varphi_{\beta}(t)&\geq&C(t_1)\omega_\beta\geq\frac{C(t_1)}{2}\omega_0.
\end{eqnarray*}
\begin{eqnarray*}
\sqrt{-1}\partial\overline{\partial}|s|_h^{2q}=q^2|s|_h^{2q}\sqrt{-1}\partial\log|s|_h^{2}\wedge\overline{\partial}\log|s|_h^{2}
+q|s|_h^{2q}\sqrt{-1}\partial\overline{\partial}\log|s|_h^{2},
\end{eqnarray*}
we obtain the estimate
\begin{eqnarray*}
\hat{\Delta}|s|_h^{2q}&\geq &q|s|_h^{2q}\int_0^1g_{sH_\beta(t)+(1-s)\varphi_\beta(t)}(\frac{\partial^2}{\partial z^i\partial\bar{z}^j}\log|s|_h^{2})ds\\
&=&-q|s|_h^{2q}\int_0^1g_{sH_\beta(t)+(1-s)\varphi_\beta(t)}^{i\bar{j}}\theta_{i\bar{j}}ds\\
&\geq&-C(t_1)q|s|_h^{2q}g_{\beta}^{i\bar{j}}g_{0,i\bar{j}}\geq-C(t_1)
\end{eqnarray*}
on $M\setminus D$, where constant $C(t_1)$ independent of $a$. Then we obtain
\begin{eqnarray*}
  \frac{\partial \Psi(t)}{\partial t}\leq\hat{\Delta}\Psi(t)-\Psi(t)+aC(t_1).
\end{eqnarray*}

Let $\tilde{\Psi}=e^{(t-t_1)}\Psi-aC(t_1)e^{(t-t_1)}-\varepsilon (t-t_1)$. By choosing suitable $0<q<1$, we can  assume that the space maximum of $\tilde{\psi}$ on $[t_1,T]\times M$ is attained away from $D$. Let $(t_0,x_0)$ be the maximum point. If $t_0>t_1$, by the maximum principle, at $(t_0,x_0)$, we have
\begin{eqnarray*}
 0\leq (\frac{\partial }{\partial t}-\hat{\Delta})\tilde{\Psi}(t)\leq -\varepsilon,
\end{eqnarray*}
which is impossible, hence $t_0=t_1$. Then for $(t,x)\in [t_1,T]\times M$, we obtain
\begin{eqnarray*}
H_\beta(t)-\varphi_\beta(t)&\leq& \|H_\beta(t_1,x)-\varphi_\beta(t_1,x)\|_{L^\infty(M)}+aC(t_1)+\varepsilon T
\end{eqnarray*}
Since $\lim\limits_{t\rightarrow0^+}\|H_{\beta}(t,z)-\psi_\beta\|_{L^{\infty}(M)}=0$ and $(\ref{2016120201})$, let $a\rightarrow0$ and then $t_1\rightarrow0^+$,
\begin{eqnarray*}
H_\beta(t)-\varphi_\beta(t)\leq \varepsilon T.
\end{eqnarray*}
It shows that $H_\beta(t)\leq\varphi_\beta(t)$ after we let $\varepsilon\rightarrow0$.
For any $(t,z)\in(0,T]\times (M\setminus D)$
\begin{eqnarray}\nonumber\varphi_{\beta}(t,z)-\psi_0(z)&\geq&te^{- t}u_{\beta}+h(t)e^{-t}+\psi_\beta-\psi_0\\
&\geq&-Ct-C(1-e^{-t})+h_{1}(t)e^{-t} ,
\end{eqnarray}
where $h_{1}(t)=n(t\log t-t)e^{t}-n\int_0^te^{s}s\log s ds$, constant $C$ independent of $\beta$. Let $\beta\rightarrow0$, we have
\begin{eqnarray}\varphi(t,z)-\psi_0(z)\geq-Ct-C(1-e^{-t})+h_{1}(t)e^{-t} ,
\end{eqnarray}
There exists $\delta_2$ such that for any $t\in[0,\delta_2]$ and $z\in M\setminus D$,
\begin{eqnarray}\label{221010}\varphi(t,z)-\psi_0(z)>-\frac{\epsilon}{2}.
\end{eqnarray}
Let $\delta=\min(\delta_1,\delta_2)$, then for any $t\in(0,\delta]$ and $z\in K$,
\begin{eqnarray}-\epsilon<\varphi(t,z)-\psi_{0}(z)<\epsilon.
\end{eqnarray}
Hence, $\varphi(t)\in C^0([0,\infty)\times (M\setminus D))$. Since $\psi_\beta$ converge to $\psi_0$ in $L^1$-sense on $M$, for sufficiently small $\beta_2$, we have
\begin{eqnarray}\int_M|\psi_{\beta_2}(z)-\psi_{0}(z)|\ \omega_0^n<\frac{\epsilon}{2}.
\end{eqnarray}
By $(\ref{2210})$, $(\ref{2016120201})$ and $(\ref{221010})$, there exists $\delta$ such that for any $t\in(0,\delta)$,
\begin{eqnarray}\int_M|\varphi(t)-\psi_{0}(z)|\ \omega_0^n<\epsilon,
\end{eqnarray}
which implies $(\ref{102})$.
\QEDB

\begin{thm}\label{225} $\omega(t)=\omega_{0}+\sqrt{-1}\partial\bar{\partial}\varphi(t)$ is a long-time solution to cusp K\"ahler-Ricci flow $(\ref{CUSPKRF})$.
\end{thm}

{\bf Proof:}\ \ We should only prove that $\omega(t)$ satisfies equation $(\ref{CUSPKRF})$ in the sense of currents on $(0,\infty)\times M$.

Let $\eta=\eta(t,x)$ be a smooth $(n-1,n-1)$-form with compact support in $(0,\infty)\times M$. Without loss of generality, we assume that its compact support included in $(\delta,T)$ ($0<\delta<T<\infty$). On $[\delta,T]\times (M\setminus D)$, $\log\frac{\omega_{\beta}^{n}(t)|s|_{h}^{2(1-\beta)}}{\omega_{0}^{n}}-\psi_\beta$, $\log\frac{\omega^{n}(t)|s|_{h}^{2}}{\omega_{0}^{n}}-\psi_0$, $\varphi_{\beta}(t)-\psi_\beta$ and $\varphi(t)-\psi_0$ are uniformly bounded. On $[\delta,T]$, we have
\begin{eqnarray}\label{201503201}\nonumber&\ &\int_{M}\frac{\partial\omega_{\beta}(t)}{\partial t}\wedge\eta=\int_{M}\sqrt{-1}\partial\bar{\partial}\frac{\partial\varphi_{\beta}(t)}{\partial t}\wedge\eta\\\nonumber
&=&\int_{M}\big(\log\frac{\omega_{\beta}^{n}(t)|s|_{h}^{2(1-\beta)}}
{\omega_{0}^{n}}-\psi_\beta-(\varphi_{\beta}(t)-\psi_\beta)+h_0\big)
  \sqrt{-1}\partial\bar{\partial}\eta\\\nonumber
&\xrightarrow{\varepsilon\rightarrow 0}&\int_{M}(\log\frac{\omega^{n}(t)|s|_{h}^{2}}{\omega_{0}^{n}}-\psi_0-(\varphi(t)-\psi_0)+h_0)\sqrt{-1}\partial\bar{\partial}\eta\\
&=&\int_{M} (-Ric (\omega(t))- \omega(t) +[D])\wedge\eta.
\end{eqnarray}
At the same time, there also holds
\begin{eqnarray}\label{201503202}\nonumber\int_M \omega_{\beta}(t)\wedge\frac{\partial\eta}{\partial t}&=&\int_M \omega_{0}\wedge\frac{\partial\eta}{\partial t}+\int_M \varphi_{\beta}(t)\sqrt{-1}\partial\bar{\partial}\frac{\partial\eta}{\partial t}\\\nonumber
&\xrightarrow{\beta\rightarrow 0}&\int_M \omega_0\wedge\frac{\partial\eta}{\partial t}+\int_M \varphi(t)\sqrt{-1}\partial\bar{\partial}\frac{\partial\eta}{\partial t}\\
&=&\int_M \omega(t)\wedge\frac{\partial\eta}{\partial t}.
\end{eqnarray}
On the other hand,
\begin{eqnarray}\label{201503203}\nonumber\frac{\partial}{\partial t}\int_M \omega_{\beta}(t)\wedge\eta
&=&\int_M \varphi_{\beta}(t)\sqrt{-1}\partial\bar{\partial}\frac{\partial\eta}{\partial t}+\int_M \omega_0\wedge\frac{\partial\eta}{\partial t}+\int_M \frac{\partial\varphi_{\beta}(t)}{\partial t} \sqrt{-1}\partial\bar{\partial}\eta\\\nonumber
&\xrightarrow{\varepsilon\rightarrow 0}&
\int_M \varphi(t)\sqrt{-1}\partial\bar{\partial}\frac{\partial\eta}{\partial t}+\int_M \omega_0\wedge\frac{\partial\eta}{\partial t}+\int_M \frac{\partial\varphi}{\partial t} \sqrt{-1}\partial\bar{\partial}\eta\\
&=&\frac{\partial}{\partial t}\int_M \omega(t)\wedge\eta.
\end{eqnarray}
Combining equality
$$\frac{\partial}{\partial t}\int_M \omega_{\beta}(t)\wedge\eta=\int_{M}\frac{\partial\omega_{\beta}(t)}{\partial t}\wedge\eta+\int_M \omega_{\beta}(t)\wedge\frac{\partial\eta}{\partial t}$$
with equalities $(\ref{201503201})$-$(\ref{201503203})$, on $[\delta,T]$, we have
\begin{eqnarray}\label{201503204}\nonumber\frac{\partial}{\partial t}\int_M \omega(t)\wedge\eta&=&\int_{M} \big(-Ric (\omega(t))- \omega(t)+[D]\big)\wedge\eta\\
&\ &+\int_M \omega(t)\wedge\frac{\partial\eta}{\partial t}.
\end{eqnarray}
Integrating form $0$ to $\infty$ on both sides,
\begin{eqnarray*}\int_{(0,\infty)\times M} \frac{\partial \omega(t)}{\partial t}\wedge\eta~dt
&=&-\int_{(0,\infty)\times M} \omega(t)\wedge\frac{\partial\eta}{\partial t}~dt=-\int_{0}^{\infty}\int_M \omega(t)\wedge\frac{\partial\eta}{\partial t}~dt\\
&=&\int_{0}^{\infty}\int_{M} \big(-Ric (\omega(t))- \omega(t) + [D]\big)\wedge\eta~dt\\
&=&\int_{(0,\infty)\times M} \big(-Ric (\omega(t))- \omega(t) +[D]\big)\wedge\eta~dt.
\end{eqnarray*}
By the arbitrariness of $\eta$, we prove that $\omega(t)$ satisfies cusp K\"ahler-Ricci flow $(\ref{CUSPKRF})$ in the sense of currents on $(0,\infty)\times M$.  \QEDB

\medskip

Now we prove the uniqueness theorem.

\begin{thm}\label{228} Let $\tilde{\varphi}(t)\in C^0\big([0,\infty)\times (M\setminus D)\big)\bigcap C^\infty\big((0,\infty)\times(M\setminus D)\big)$ be a long-time solutions to parabolic Monge-Amp\`ere equation
\begin{eqnarray}\label{2016120301}
  \frac{\partial \varphi(t)}{\partial t}=\log\frac{(\omega_{0}+\sqrt{-1}\partial\bar{\partial}\varphi(t))^{n}}{\omega_{0}^{n}}-\varphi(t)+h_{0}
  +\log|s|_{h}^{2}
\end{eqnarray}
on $(0,\infty)\times (M\setminus D)$. If $\tilde{\varphi}$ satisfies
\begin{itemize}
  \item For any $0<\delta<T<\infty$, there exists uniform constant $C$ such that
\begin{eqnarray*}C^{-1}\omega_{cusp}\leq \omega_{0}+\sqrt{-1}\partial\bar{\partial}\tilde{\varphi}(t)\leq C\omega_{cusp}\ \ \ on\ \ [\delta,T]\times (M\setminus D);
\end{eqnarray*}
\item on $(0,T]$, $\|\tilde{\varphi}(t)-\psi_0\|_{L^{\infty}(M\setminus D)}\leqslant C$;\\
  \item on $[\delta, T]$, there exist constant $C^{*}$ such that $\| \frac{\partial\tilde{\varphi}(t)}{\partial t}\|_{L^{\infty}(M\setminus D)}\leqslant C^{*}$;
  \item $\lim\limits_{t\rightarrow0^{+}}\|\tilde{\varphi}(t)-\psi_{0}\|_{L^{1}(M)}=0$.
  \end{itemize}
Then $\tilde{\varphi}(t)\leq\varphi(t)$.
\end{thm}

{\bf Proof:}\ \  For any $0<t_1<T<\infty$ and $a>0$. Denote $\Psi(t)=\tilde{\varphi}(t)+a\log|s|_h^2-\varphi_\beta(t)$ and $\hat{\Delta}=\int_0^1 g_{s\tilde{\varphi}(t)+(1-s)\varphi_\beta(t)}^{i\bar{j}}\frac{\partial^2}{\partial z^i\partial\bar{z}^j}ds$. We note that $\tilde{\varphi}(t)$ is bounded from above. $\Psi(t)$ evolves along the following equation
\begin{eqnarray*}
  \frac{\partial \Psi(t)}{\partial t}=\hat{\Delta}\Psi(t)-a\hat{\Delta}\log|s|_h^{2}-\Psi(t)+(a+\beta)\log|s|_h^2.
\end{eqnarray*}
Since $-\sqrt{-1}\partial\bar{\partial}\log|s|_h^2=\theta$, we obtain
\begin{eqnarray*}
-\hat{\Delta}\log|s|_h^2=\int_0^1g_{s\tilde{\varphi}(t)+(1-s)\varphi_\beta(t)}^{i\bar{j}}\theta_{i\bar{j}}ds\leq C(t_1)
\end{eqnarray*}
on $M\setminus D$. Then we obtain
\begin{eqnarray*}
  \frac{\partial \Psi(t)}{\partial t}\leq\hat{\Delta}\Psi(t)-\Psi(t)+aC(t_1).
\end{eqnarray*}
Then by the arguments as that in Proposition $\ref{101}$, on $[t_1,T]\times (M\setminus D)$,
\begin{eqnarray*}
\tilde{\varphi}(t)-\varphi_\beta(t)\leq e^{-(t-t_1)}\sup_M(\tilde{\varphi}(t_1)-\varphi_\beta(t_1)).
\end{eqnarray*}
Since $\tilde{\varphi}(t_1)$ converge to $\psi_0$ in $L^1$-sense and $\varphi_\beta(t)$ converge to $\psi_\beta$ in $L^\infty$-sense as $t_{1}\rightarrow0^+$, by Hartogs Lemma, we have
\begin{eqnarray*}
\tilde{\varphi}(t)-\varphi_\beta(t)\leq e^{-t}\sup_M(\psi_0-\psi_\beta)\leq0,
\end{eqnarray*}
after we let $t_1\rightarrow0$. Hence $\tilde{\varphi}(t)\leq\varphi(t)$ on $(0,\infty)\times(M\setminus D)$.\QEDB

\begin{rem} If $M$ is a compact K\"ahler manifold with smooth hypersurface $D$. We can also consider unnormalized cusp K\"ahler-Ricci flow
\begin{eqnarray}\label{NCUSPKRF}
\begin{cases}
  \frac{\partial \hat{\omega}(t)}{\partial t}=-Ric(\hat{\omega}(t))+[D].\\
  \\
  \hat{\omega}(t)|_{t=0}=\omega_{cusp}\\
  \end{cases}
  \end{eqnarray}
If we define $\omega(t)=e^{-t}\hat{\omega}(e^t-1)$, then flow $(\ref{NCUSPKRF})$ is actually the same as normalized cusp K\"ahler-Ricci flow $(\ref{CUSPKRF})$ only moduli a scaling. Let
\begin{eqnarray}T_0=\sup\{\ t\ |\ [\omega_0]-t(c_1(M)-c_1(D))>0\}.
\end{eqnarray}
Combining the arguments of Tian-Zhang \cite{GTZZ} and Liu-Zhang \cite{JWLXZ1} with the arguments in this paper, there exists a unique solution to flow $(\ref{NCUSPKRF})$ on $[0,T_0)$ in some weak sense which is similar as Definition \ref{04.5}.
\end{rem}

\section{The convergence of cusp K\"ahler-Ricci flow}
\setcounter{equation}{0}

In this section, we prove the convergence theorem of cusp K\"ahler-Ricci flow $(\ref{CUSPKRF})$.

{\bf Proof of Theorem \ref{05}:}\ \ Differentiating equation $(\ref{2016120301})$ in time $t$, we have
\begin{eqnarray}
(\frac{d}{dt}-\Delta_t)\frac{\partial\varphi}{\partial t}=-\frac{\partial\varphi}{\partial t}
\end{eqnarray}
on $[\delta,T]\times (M\setminus D)$ with $\delta>0$. For any $\varepsilon>0$,
\begin{eqnarray}\nonumber
(\frac{d}{dt}-\Delta_t)(\frac{\partial\varphi}{\partial t}+\varepsilon\log|s|_h^2)&=&-\frac{\partial\varphi}{\partial t}+\varepsilon tr_{\omega(t)}\theta\\
&\leq&-(\frac{\partial\varphi}{\partial t}+\varepsilon\log|s|_h^2)+\varepsilon C(\delta,T),
\end{eqnarray}
where constant $C(\delta,T)$ independent of $\varepsilon$. For any $\eta>0$, let $H=e^{t-\delta}(\frac{\partial\varphi}{\partial t}+\varepsilon\log|s|_h^2)-\varepsilon e^{t-\delta} C(\delta,T)-\eta(t-\delta)$. Since $\frac{\partial\varphi}{\partial t}$ is bounded on
$[\delta,T]\times (M\setminus D)$, the maximum point $(t_0,x_0)$ of $H$ satisfies $x_0\in M\setminus D$. By the maximum principle, $t_0=\delta$. Hence,
\begin{eqnarray}
\frac{\partial\varphi}{\partial t}\leq C(\delta)e^{-t}-\varepsilon\log|s|_h^2+\varepsilon C(\delta,T)+\eta T.
\end{eqnarray}
Let $\varepsilon\rightarrow0$, $\eta\rightarrow0$ and then $T\rightarrow\infty$, we obtain
\begin{eqnarray}
\frac{\partial\varphi}{\partial t}\leq C(\delta)e^{-t}\ \ \ on\ \ [\delta,\infty)\times (M\setminus D).
\end{eqnarray}
By the same arguments, we can get the lower bound of $\frac{\partial\varphi}{\partial t}$. In fact, we obtain
\begin{eqnarray}
|\frac{\partial\varphi}{\partial t}|\leq C(\delta)e^{-t}\ \ \ on\ \ [\delta,\infty)\times (M\setminus D).
\end{eqnarray}
For $\delta<t<s$,
\begin{eqnarray}
|\varphi(t)-\varphi(s)|\leq C(\delta)(e^{-t}-e^{-s})\ \ \ on\ \ [\delta,\infty)\times (M\setminus D).
\end{eqnarray}
Therefore, $\varphi(t)$ converge exponentially fast in $L^\infty$-sense to $\varphi_\infty$ on $M\setminus D$. And $\varphi(t)$ converge to $\varphi_\infty$ in $C^\infty_{loc}$-sense on $M\setminus D$. At the same time, For any smooth $(n-1,n-1)$-form $\eta$,
\begin{eqnarray}
\int_M\frac{\partial\omega(t)}{\partial t}\wedge\eta=\int_M\frac{\partial\varphi(t)}{\partial t}\sqrt{-1}\partial\bar{\partial}\eta\xrightarrow{t\rightarrow \infty}0
\end{eqnarray}
while
\begin{eqnarray*}
\int_M\frac{\partial\omega(t)}{\partial t}\wedge\eta&=&\int_M\sqrt{-1}\partial\bar{\partial}(\log\frac{|s|_h^2(\omega_{0}+\sqrt{-1}\partial\bar{\partial}\varphi(t))^{n}}{\omega_{0}^{n}}-\varphi(t)+h_{0})\wedge\eta\\
&=&\int_M(\log\frac{|s|_h^2(\omega_{0}+\sqrt{-1}\partial\bar{\partial}\varphi(t))^{n}}{\omega_{0}^{n}}-\psi_0-(\varphi(t)-\psi_0)+h_{0})\sqrt{-1}\partial\bar{\partial}\eta\\
&\xrightarrow{t\rightarrow \infty}&\int_M(\log\frac{|s|_h^2(\omega_{0}+\sqrt{-1}\partial\bar{\partial}\varphi_\infty)^{n}}{\omega_{0}^{n}}-\psi_0-(\varphi_\infty-\psi_0)+h_{0})\sqrt{-1}\partial\bar{\partial}\eta\\
&=&\int_M(-Ric(\omega_\infty)-\omega_\infty+[D])\wedge\eta.
\end{eqnarray*}
which implies the convergence in the sense of currents. \QEDB

\hspace{1.4cm}

\end{document}